# The Use of Harmonic Information in the Optimal Synthesis of Mechanisms

A.M.CONNOR[1], S.S.DOUGLAS & M.J.GILMARTIN

**SUMMARY** *This paper reviews several uses of harmonic information in the synthesis of mechanisms and shows that such information can be put to even greater use in this field. Results are presented for both single and multi-degree of freedom systems which support this claim. In both cases, the inclusion of harmonic information into the objective function aids the search to locate high quality solutions.*

## 1. Introduction

In the synthesis of mechanisms the accepted methods is to apply kinematic optimisation followed by possibly dynamic optimisation. In this approach the kinematic design requirements to be optimised are usually structural error and transmission angle. That is the emphasis is placed on satisfying the desired motion. A major drawback to this approach is that the resultant geometry may lead to solutions that are dynamically unattractive This paper proposes an alternative to existing methods by also including harmonic information into the optimisation criteria. The reason for this is that by through the study of the harmonic content an insight into the dynamic performance of mechanisms can be obtained. For example compare the case were two mechanisms have been synthesised for some design requirement. One solution possess a substantial harmonic content, over the other which mechanisms should be adopted? If the mechanism with the substatial harmonic content is chosen it is reasonable to expect that these harmonics would excite unwanted vibrations. If the other mechanism was chosen then the potential for adverse vibrations to be stimulated is much reduced and a superior dynamic performance obtained. This paper demonstrates how single and multi degree of freedom mechanisms can be synthesised to exhibit good motions via the optimisation of harmonic spectra.

## 2. Literature Review

Much of the previously published work which utilises harmonic information has been concerned with the analysis of the output motions in terms of harmonic series [3]. Another example of early work [4] developed dynamic performance criteria for crank rocker linkages. The use of harmonic analysis revealed that the higher harmonics of the accelerations were considerably smaller than the fundamental harmonic. In all cases, the method found results where the ratio of the input crank to the ground link was small.

Some work has included such information in the synthesis of mechanisms. Some examples of this [5,6,7] make use of the observation that mechanisms with short input cranks are likely to produce output motions that are relatively small. These motions can therefore be described as approximations of simple harmonic functions of the input crank angle. This work has also been extended to cover multi-degree of freedom mechanisms[8].

Other approaches have been developed which include harmonic information as a basis for mechanism synthesis. One example [9] of this is the use of harmonic coefficients to store mechanism information in a library. Comparing the harmonic coefficients of the desired motion to those in the library enables a mechanism configuration which approximates the required motion to be selected. This can then be numerically optimised to reduce the error between the actual and the desired output.

## 3. The Fourier Transform

---

[1]    A.M.Connor, Engineering Design Centre, School of Mechanical Engineering, University of Bath, Bath, BA2 7AY.
       +44 (0)1225 826131, ensamc@bath.ac.uk, http://www2.bath.ac.uk/~ensamc/

Whilst the results presented in later sections are for different applications, and harmonic information is used in different ways, both methods utilise a discrete Fourier transform (DFT). Essentially, the DFT represents time domain information as a complex series. The DFT is derived from the analytical Fourier transform. The only difference between the DFT and analytical techniques is that the coefficients relating to the Fourier series are calculated by numerical integration.

The required series is denoted by;

$$f(x) = \frac{1}{2}a_0 + \sum_{n=1}^{\infty}\{a_n \cos nx + b_n \sin nx\}$$

The individual coefficients are calculated using the following expressions.

$$a_0 = \frac{1}{\pi}\int_0^{2\pi} f(x)dx$$

$$a_n = \frac{1}{\pi}\int_0^{2\pi} f(x)\cos nx\, dx$$

$$b_n = \frac{1}{\pi}\int_0^{2\pi} f(x)\sin nx\, dx$$

Each harmonic, other than the fundamental harmonic, consists of a real and imaginary part. The magnitude of each harmonic is given as;

$$H_n = \sqrt{a_n^2 + b_n^2}$$

Similarly, the phase of each harmonic is given as;

$$\phi(H_n) = \tan^{-1}\left(\frac{b_n}{a_n}\right)$$

The importance of the magnitude and phase of each harmonic can be realised by considering a hypothetical series. Each term in the series represents a sinusoid in the time domain. Adding the total number of sinusoids, taking into account phase difference, using the principle of superposition results in the time domain function.

**4. Single Degree of Freedom Systems**

The origins of this work was an experimental study of industrial sewing mechanisms. This revealed that the mechanisms significant harmonic content was causing unacceptable stitch patterns at high speeds.

The results presented demonstrate how harmonic synthesis ideas have been applied to single degree of freedom industrial sewing mechanisms. In both cases, the optimisation method was a hill climbing technique. using Powells method of conjugate directions [10] and a scaled lagrangian penalty function [11].

4.1 *Straight Line Needle Mechanism*

This mechanism consists of a four bar crank rocker coupled to a slider crank. The crank rotation provides the input motion and the slider oscillation provides the output motion. The kinematics objective of this design was to determine link geometry such that a desired translation of the needle is produced for a fully rotational crank input. To this objective the design requirement that the solution possess smooth motion and good transmission angle characteristics was added. This became the optimisation criteria. The output translation requirements became the equality constraints and the mechanism feasibility was monitored by the inequality constraints. Figure 1 shows the trial mechanism geometry.

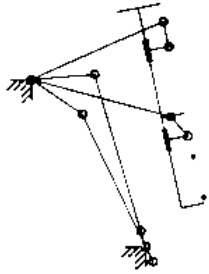

Figure 1 : Trial Mechanism

Nine design parameters were selected to be freed for the optimisation. Figure 2 displays the two limit positions of the optimised solution. It can be seen that the mechanism strokes between the specified equality constraints (shown as diamonds) and also that good transmission angles have been determined.

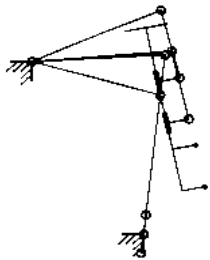

Figure 2 : Optimised Solution

The table in Figure 3 shows the harmonic contributions to the output motion. It can be seen that motion consists principley of the fundamental harmonic. This clearly indicates that the optimisation criteria which penalises the square of the second to fifth harmonics has been successful.

| Order | Magnitude | Phase Angle |
|---|---|---|
| 1 | 12.2546 | 186.782 |
| 2 | 0.3852 | 102.837 |
| 3 | 0.3181E-01 | 295.258 |
| 4 | 0.2136E-02 | 227.743 |
| 5 | 0.5629E-04 | 87.6564 |

Figure 3 : Harmonic Information

### 4.2  Right Hand Looper Mechanism

In this synthesis the objective was to design a spatial - planar mechanism. The spatial (RRRSR) mechanism converted a continuously rotating crank input into an orthogonal output oscillation. The planar mechanism design requirement was to generate a coupler curve.

The synthesis problem was partitioned into two related optimisation strategies. In the first the kinematics objective of the design was to determine a four bar coupler curve mechanism which gave good approximation to a specified segment of a coupler curve. For this design optimisation criterion comprising coupler curve structural error and smooth output motion and good transmission angles was selected. As the objective function was multi functional each criteria was weighted to enable a range of searches to be undertaken. In this synthesis formulation the input oscillation is unknown and inequality constraints were constructed to keep the mechanism feasible over its range of operation. Once an optimal solution has been found the derived input oscillation became the equality constraints for the spatial mechanism optimisation. The objective criteria of the spatial drive was again for smooth motion on the output oscillating link and also good spatial transmission angle on the output dyad. The inequality constraints have been constructed to keep the solution feasible.

Figure 4 shows a trial mechanism and the boundary defining the feasible zone for the coupler curve. The vertical dashed line represents the output path of the first optimisation and the diamonds represent the desired coupler curve.

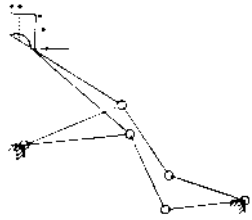

Figure 4 : Trial Mechanism

In this synthesis because the input oscillation is unknown an iterative technique is applied to determine the input angular throw which positions the trial mechanism as close as possible to the desired coupler curve. Figure 5 shows the optimal solution obtained when six mechanism parameters were freed. In this optimisation the relative weighting between the structural error and the harmonic and transmission functions was one hundred to one.

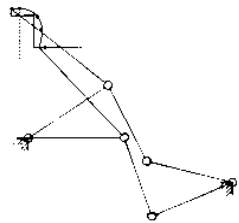

Figure 5 : Optimised Solution

Figure 6 shows the optimised spatial drive. In this optimisation the spatial transmission angle and the smooth harmonic criteria were equally weighted.

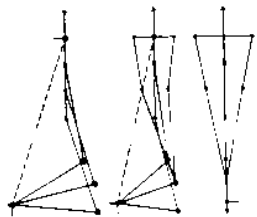

Figure 6 : Optimised Spatial Mechanism

**5. Multi-Degree of Freedom Systems**

The results presented in this section illustrate how harmonic information can be used in the synthesis of multi-degree of freedom mechanisms. The particular application is that of a hybrid five bar path generating mechanism [12]. Hybrid mechanisms have been defined [13,14,15] as a multi-degree of freedom mechanism, where different motor types provide the inputs required. All research in this area to date concerns mechanisms with two degrees of freedom, where one input is provided by a constant velocity (CV) motor and the other input is provided by a programmable servo motor. The desired optimal condition for such configurations is considered to be when the CV motor provides the bulk of the motion. When this occurs, the size of the servo motor is small and there is large potential for energy regeneration on the CV motor axis by the inclusion of a flywheel.

Previous work [15] has utilised objective function criteria based on either torque requirement, or power consumption, and distribution. Such objective functions are computationally expensive. This work attempts to replace such dynamic criteria with an approximation which utilises frequency domain information.

Minimising the magnitude of high order (higher frequency) harmonics of the servo motor profile is likely to lead to high quality solutions for the following reasons. As the magnitude of low orders is reduced, the motion profile trends towards simple harmonic motion. Simple harmonic motion is infinitely differentiable without introducing discontinuities in higher derivatives. Therefore, a motor profile with low order harmonics is likely to exhibit desirable velocity, acceleration and torque profiles.

The search method used in this instance is based around a Genetic Algorithm [17]. GAs are a novel search and optimisation method based on natural selection, population genetics and the principle of survival of the fittest. This method has been shown to be readily applicable to mechanism synthesis problems [2,12,18,19].

5.1 *Five Bar Mechanisms*

A fivebar mechanism is a two degree of freedom mechanism which requires two inputs to fully define the mechanism motion. Figure 7 illustrates such a mechanism.

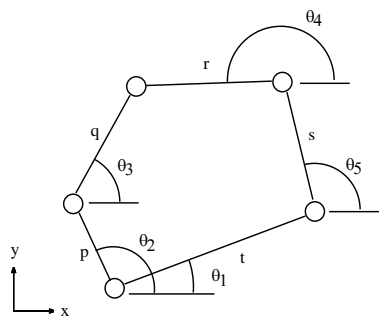

Figure 7 : Five Bar Mechanism

In a hybrid mechanism of this type, one input is provided by a programmable servo motor whilst the other is provided by a CV motor. In Figure 7, $\theta_1$ is considered to be the CV motor input whilst $\theta_5$ is the servo motor input. Link t is the ground link and all motions are expressed relative to the x-axis..

However, if the synthesis is to search for the mechanism link lengths which result in servo motor profile being the most desirable, then obviously the motor motion profile is not known. It is therefore necessary to specify alternative inputs for use in the analysis routine embedded in the synthesis method. This is achieved by considering the nature of the mechanism and relating this to the desired output motion of the mechanism. For the sake of the analysis and synthesis of the five bar mechanism, it is assumed that the end effector is the revolute joint between links q and r. The input motion for link p is known, and so by defining the actual position of the end effector it is possible to fully describe the motions in the mechanism.

5.2 *Output Motion Error*

Figure 8 illustrates how the positions of links p and q can be used to define the required angles of the mechanism so that it may be analysed fully and also calculate an error score.

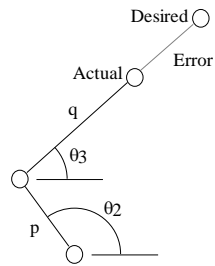

Figure 8 : Calculation of Error

A vector is defined from the end of link p to the desired position of the end effector for the given value of $\theta_2$. The actual position of this point is then calculated, given the length of link q, and the error between the two points found. This error is summed around the cycle for the input link and forms the basis of one component of the objective function.

$$obj_{err} = error^2$$

5.3 *Mechanism Mobility*

The two links, p and q, form a dyad which can be used to evaluate the error at a given point for a given input angle. Similarly, the remaining two mobile links, r and s, form another dyad which is used to calculate a penalty function criteria for use in the overall objective function.

This penalty function is based on the mobility of the mechanism. For a truly mobile mechanism, the dyad formed by the links r and s should be able to 'close' for all given positions of the CV input crank. This means that the position of the common revolute joint, when considered as part of the dyad formed by r and s, should be able to reach the actual position defined by the dyad formed by the links p and q.

For each position that this dyad cannot close, the mobility counter is increased. It is important to realise that for each position of the input link, there are two possible closures of the dyad. This is illustrated in Figure 9.

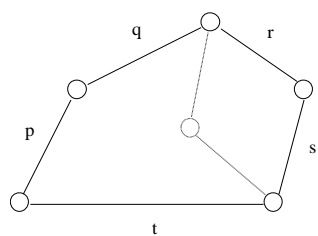

Figure 9 : Multiple Closures of the Mechanism

The implication of this, at this stage, is that for each position of the CV input link the mobility score may increase by two. The mobility penalty score is summed throughout the cycle and, for twenty four precision points, may range between zero and forty eight. As with the error score value, it has been found that raising the mobility penalty to a power acts so that the search is directed more towards feasible solutions. The following expression is used to calculate the contribution of the mobility penalty function to the overall objective function.

$$obj_{mob} = mobility^3$$

Once two inputs to the mechanism are known, in this case $\theta_2$ and $\theta_3$, it is possible to calculate any other angle in the mechanism. In this case, $\theta_5$ is calculated using the equations derived in section 4.3.

The values for $\theta_5$ around the cycle of the mechanism can be used to calculate a number of figures of merit which describe the quality of the motion profile. However, the situation is complicated by the presence of multiple closures in the mechanism. For each input position, two values for $\theta_5$ can be calculated and the question which arises is "which value should be chosen?".

An algorithm has been developed to choose the set of servo motor displacements for a given mechanism. The algorithm is based on the fact that it is desirable to have smooth velocity profiles for the servo motor input.

### 5.4 Closure Tracking Algorithm

The algorithm minimises the fluctuation in servo motor velocity in the following way. Given that the mechanism is assumed to start in the 'open' closure, the algorithm chooses the next displacement angle (from a choice of two) based on the magnitude of the difference in displacement between each possible angle and the starting angle. For each subsequent point, the algorithm chooses a displacement by considering the difference in magnitude and direction between each possible velocity for the current stage and known velocity of the previous stage. The effect of this is that the at each step of the cycle, the closure is chosen so that the trend in displacements continues with the smallest change in displacement. The trend is only broken if the angle in the next step produces a step change many times larger than the alternative which results in a change of direction.

### 5.5 Assessment of Quality

Given that the servo motor displacement is now known, it is possible to calculate a figure of merit which describes the quality of the motion profile. Previous work [2] has shown that the use of only harmonic content does not force the search towards dynamically optimal mechanisms. This is because the method locates a mechanism where the closing dyad (link r and s) has long links when compared to the CV input crank. This occurs because such long links produce servo motor motion profiles where the a large cartesian displacement is produced by a small angular displacement. Because of this the fundamental harmonic, which describes the magnitude of the displacements in the time domain, becomes significant when compared to the values of low order harmonics. To overcome this, a second figure of merit has been developed which works in conjunction with the harmonic content. These two design objectives can now be described in more detail.

#### 5.5.1 Harmonic Content

One of the design objectives used is based on the magnitude of the harmonics of the motor displacement profile. The contribution to the overall objective function is calculated by summing the value of each harmonic raised the power of its order plus one. This effectively penalises low order harmonics with high magnitudes.

$$obj_{harm} = \sum_{i=1}^{n} \left( \sqrt{a_n^2 + b_n^2} \right)^{n+1}$$

#### 5.5.2 Motion Swept Area

Using only the harmonic content of the motion profile as a design objective does not lead to approximate dynamically optimum mechanisms. This is because the method is prone to finding solutions that have long links in the closing dyad which require greater torque to move them.

Whilst short input links are desirable, it is also important to take into account the motion profile in time domain form. This is achieved by calculating the swept area of the servo motor input. The swept area is defined as being the RMS of the motor displacement multiplied by the length of link s.

$$obj_{swept} = s \sqrt{\frac{\sum_{i=1}^{n} \theta_{5i}^2}{n}}$$

5.6  *Problem Representation*

The required motion is defined by twenty four discrete points. The curve defined by this points exhibits two cusps, or points of zero velocity, and so may be for an application such as a pick and place mechanism. The parameters used in the search for the optimal mechanism are the lengths of links p,q,r and s. The length of link t is defined by including in the search the co-ordinates of the ground inputs. These are defined as local (x,y) co-ordinates with in a global constraint envelope defined by the user.

The search was terminated when the structural error dropped below a given threshold, but only if the link lengths of the located mechanism were such that no link was more than five times longer than the length of the CV motor input crank.

5.7  *Results*

Five trial runs of the problem were carried out. In each case a very rapid initial reduction in objective function value was followed by a slow, but definite, improvement.

Each of the solutions found was different due to the probabilistic nature of the GA used as the search algorithm. However, each solution had similar objective function values. The best solution was obtained in the fifth run and had the following parameter set.

$(x,y)_{cv}$ = 0,1
$(x,y)_{servo}$ = 27,-2
p = 16
q = 24
r = 30
s = 16

The graph in Figure 10 shows the required servo motor input motion.

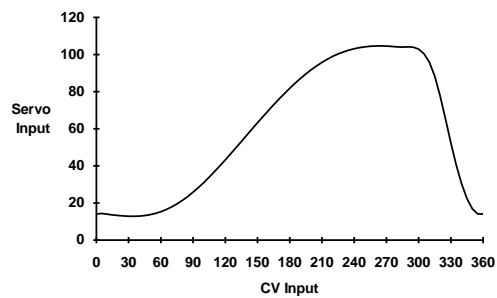

Figure 10 : Servo Motor Displacement

Figure 11 shows the actual output motion of the end effector in relation to the desired motion.

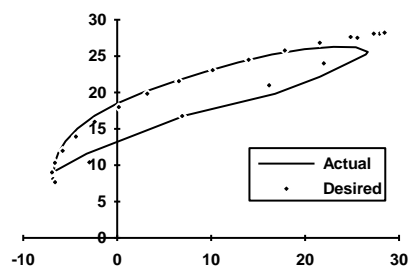

Figure 11 : End Effector Motion

The graph in Figure 12 shows the torque requirement for the mechanism. The mechanism can be seen to be operating as a valid hybrid device as the torque requirement of the servo motor is lower than that for the CV motor.

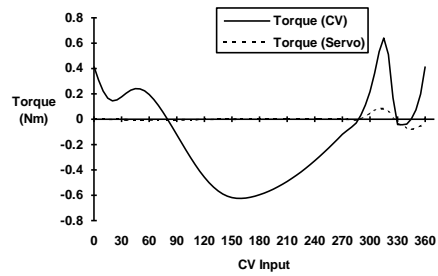

Figure 12 : Mechanism Torque Requirements

## 6. Discussion

The results presented in the previous sections show that the use of harmonic information in the synthesis of mechanisms can force the search method to find high quality solutions.

In the single degree of freedom examples, reducing the harmonic content of the output motion allowed higher operating speeds to be achieved.

In the multi-degree of freedom example, introducing harmonic information has eliminated the necessity for computationally expensive dynamic optimisation criteria and has produced a mechanism which exhibits a very desirable torque distribution between the two axes of control.

## 7. Conclusions

The results presented show the usefulness of harmonic information in the synthesis of mechanisms. This information has been applied in two examples to show how improved dynamic performance can be achieved.